\documentclass[12pt]{amsart} 
\usepackage{amssymb,amsfonts}
\usepackage[all,arc]{xy}
\usepackage{enumerate}
\usepackage{mathrsfs}
\usepackage[linktoc=all]{hyperref}
\hypersetup{
	colorlinks,
	citecolor=black,
	filecolor=black,
	linkcolor=black,
	urlcolor=black
}

\newtheorem{thm}{Theorem}[section]

\newtheorem{conj}[thm]{Conjecture}

\theoremstyle{definition}
\newtheorem{defn}[thm]{Definition}

\theoremstyle{remark}
\newtheorem{rem}[thm]{Remark}

\makeatletter
\let\c@equation\c@thm
\makeatother
\numberwithin{equation}{section}

\usepackage{cite}
\title{Area, Scalar Curvature, and Hyperbolic 3-manifolds} 
\author{Ben Lowe}
\email{benl@princeton.edu} 

\begin{document}
	\maketitle 
	\begin{abstract} 
		Let $M$ be a closed hyperbolic 3-manifold that admits no infinitesimal conformally-flat deformations.  Examples of such manifolds were constructed by Kapovich.   Then if $g$ is a Riemannian metric on $M$ with scalar curvature greater than or equal to $-6$, we find lower bounds for the areas of stable immersed minimal surfaces $\Sigma$ in $M$.  Our bounds improve the closer $\Sigma$ is to being homotopic to a totally geodesic surface in the hyperbolic metric.  
		
		We also consider a functional introduced by Calegari-Marques-Neves
	    that is defined by an asymptotic count of minimal surfaces in $(M,g)$. We show this functional to be uniquely maximized, over all metrics of scalar curvature greater than or equal to $-6$, by the hyperbolic metric.    Our proofs use the Ricci flow with surgery.  
		

		

		
		
		
		\end{abstract}

	\section{Introduction} 
	
	\subsection{} 
	
	There are many ``infinitesimal-to-global" theorems in geometry that bound metric quantities associated to a Riemannian manifold under assumptions on curvature.  A Riemannian manifold that realizes equality for such a bound must often be very symmetric.   Our goal here will be to obtain bounds on the areas of certain surfaces assuming a negative lower bound on the scalar curvature of the ambient metric. We will show that equality in some of the inequalities we obtain implies that the metric had to have constant sectional curvature -1.

	It follows from the second variation formula, the Schoen-Yau rearrangement trick, and the Gauss-Bonnet formula that a stable minimal surface in a 3-manifold with scalar curvature greater than or equal to $-6$ has area at least $-2\pi\chi(\Sigma)/3$. Nunes  proved that equality implies that the metric  splits as a Riemannian product of a constant negative curvature metric on a surface and an interval \cite{nunes}.    
	In this paper we improve these results for surfaces in closed hyperbolic 3-manifolds that admit no infinitesimal conformally-flat deformations.

	\subsection{} 
Let $M$ be a closed hyperbolic 3-manifold and let $\Sigma$ be a closed immersed surface in $M$.   Let $\mathcal{F}_\Sigma$ denote the set of all immersed surfaces $\Sigma'$ in $M$ homotopic to the inclusion of $\Sigma$ in $M$--- i.e., for which there exists a map $F:\Sigma \times [0,1] \rightarrow M$ such that $F|_{\Sigma \times \{0\})}$ is the inclusion of $\Sigma$ in $M$ and $F(\Sigma \times \{1\})=\Sigma'$.  For a Riemannian metric $g$ on $M$, define
\[
\mathcal{A}_{\Sigma}(M,g) = \inf \{ \text{area}(\Sigma,g) : \Sigma \in \mathcal{F}_\Sigma\}.  
\]

	 	It follows from \cite{km}, \cite{kmcounting} and \cite{seppi}  that every closed hyperbolic 3-manifold contains lots of closed immersed $\pi_1$-injective minimal surfaces whose principal curvatures are pointwise as small as desired.  Fix such a surface $\Sigma$, and set 
	 \[
	 	\epsilon_{\Sigma} = \frac{1}{\text{Area}_{g_{hyp}}(\Sigma)} \int_{\Sigma} |A|^2 ,
	 	\] 
	 	where $A$ is the second fundamental form of $\Sigma$.    The previous version of this paper contained an incorrect proof of the following statement, which we now state as a conjecture.  
	 
	 	\begin{conj} \label{varcurvintroconj} 
	 		Suppose that the scalar curvature of $g$ is greater than or equal to $-6$.  Then for every $\delta>0$ there exists $\epsilon$ (depending on $g$) such that if $\epsilon_{\Sigma}<\epsilon$, then
	 		\begin{equation} \label{varcurvineq} 
	 		\mathcal{A}_{\Sigma}(M,g) \geq (1-\delta)   \mathcal{A}_{\Sigma}(M,g_{hyp}).      
	 		\end{equation} 
	 		
	 	\end{conj} 
 	
 	We say that $(M,g_{hyp})$ admits no \textit{infinitesimal conformal deformations} if the cohomology group 
 	\begin{equation} \label{tangspace} 
 	H^1(\pi_1(M),\text{Ad}) 
 	\end{equation} 
 	\noindent vanishes, where $\text{Ad}$ is the Adjoint representation of $\pi_1(M)\subset SO(3,1)$ on $so(4,1)$ via the inclusion $so(3,1) \hookrightarrow so(4,1)$.   The vanishing of (\ref{tangspace}) implies but is strictly stronger than (see \cite{sca02}) the statement that there are no non-trival deformations of $g_{hyp}$ through conformally flat metrics.  Kapovich gave infinitely many examples of closed hyperbolic 3-manifolds for which (\ref{tangspace}) vanishes \cite{kap94}. These are obtained by Dehn surgeries on hyperbolic 2-bridge knots.  (A knot is 2-bridge if it can be realized by an embedding in $\mathbb{R}^3$ with only two local maxima.  The figure-eight knot is an example)

 	\begin{thm} \label{varcurvintro} 
 		Conjecture \ref{varcurvintroconj} is true if $(M,g_{hyp})$ admits no infinitesimal conformal deformations.  
 		\end{thm} 

	 	\subsection{}

	  Calegari-Marques-Neves \cite{cmn} recently introduced a functional on Riemannian metrics on $M$ based on an asymptotic count of minimal surfaces homotopic to immersed almost totally geodesic surfaces in the hyperbolic metric, and showed that the hyperbolic metric uniquely minimizes this functional over all metrics with sectional curvature at most $-1$.

	  To state their result, we need to introduce some terminology.  Let $S(M)$ denote the set of subgroups of $\pi_1(M)$ isomorphic to the fundamental group of a closed surface, or \textit{surface subgroups}, up to the equivalence relation of conjugacy in $\pi_1(M)$. The \textit{limit set} of a surface subgroup $\pi_1(\Sigma)$ of $\pi_1(M)$ is the set of accumulation points in $\partial_\infty \mathbb{H}^3$ of any orbit of the action of $\pi_1(\Sigma)$ on $\mathbb{H}^3$.

	  A homeomorphism $f:S^2 \rightarrow S^2$ is \textit{K-quasiconformal} if for any ball $B(x,r)\subset S^2$ there exists $r'>0$ such that $B(f(x),r')\subset f(B(x,r)) \subset B(f(x),Kr')$.  Here $B(x,s)$ denotes the ball centered at $x$ of radius $s$ in the round metric on $S^2$.  A \textit{K-quasicircle} is a Jordan curve that is the image of a round circle under a K-quasiconformal map $S^2 \rightarrow S^2$.  Kahn-Markovic showed the existence of surface subgroups of $\pi_1(M)$ whose limit sets in $\partial_\infty \mathbb{H}^3 \cong S^2$ are $K$-quasicircles for $K \searrow 1$.  It will be important for us that a surface subgroup of $\pi_1(M)$ whose limit set is a $(1+\epsilon)$-quasicircle is realized by the injective inclusion of a minimal surface with principal curvatures bounded above by some universal constant times $\log(1+\epsilon)$ \cite{seppi}.  
	  
	  Denote by $S_{\epsilon}(M)$ the elements of $S(M)$ whose limit sets are $1+\epsilon$-quasicircles. Elements of $S(M)$ are in one-to-one correspondence with homotopy classes of $\pi_1$-injective immersed surfaces in $M$.  For a metric $g$ on $M$ and for $S \in S(M)$, we denote by $\text{Area}_g(S)$ the infimal area with respect to $g$ of a surface in the corresponding homotopy class.  Calegari-Marques-Neves defined the following ``entropy" functional on metrics $g$: 
	   \[
	   E(g):= \lim_{\epsilon \to 0} \lim_{L \to \infty} \inf \frac{ \log (\# \text{Area}_g (S) \leq 4 \pi(L-1) : S \in S_{\epsilon}(M) )}{L \log L}.   
	   \]
	  
	  \noindent 
	  Their main result is the following: 
	  
	  \begin{thm} [Calegari-Marques-Neves] \label{maincmn} 
	  	If the sectional curvature of $g$ is less than or equal to $-1$, then 
	  	\[
	  	E(g) \geq E(g_{hyp})= 2
	  	\]
	  	with equality if and only if $g$ is isometric to the hyperbolic metric $g_{hyp}$.  	
	   \end{thm}

	   \noindent Using Theorem \ref{varcurvintro}, we can prove an analogue of their theorem for scalar curvature, but with a curvature bound in the opposite direction.  The first version of this paper contained an incorrect proof of the following statement, which we now state as a conjecture.  
	  
	  \begin{conj} \label{cmnscalintroconj} 
	  	If the scalar curvature of $g$ is greater than or equal to $-6$, then 
	  	\[
	  		E(g) \leq E(g_{hyp}).  
	  	\]
	  
	 \noindent with equality if and only if $g$ is isometric to $g_{hyp}$.  	
	  \end{conj} 
  
  As above we are are able to prove it in a special case.  
  
  \begin{thm} \label{cmnscalintro}
  	Conjecture \ref{cmnscalintroconj} is true if $(M,g_{hyp})$ has no infinitesimal conformally flat deformations.  
  \end{thm}


 \subsection{} 
In our proofs, we study how areas of surfaces evolve under Ricci flow.  We use the fact, due to Perelman (\cite{perelman1},\cite{perelman2},\cite{perelman3}), that Ricci flow with surgery asymptotically takes any metric on $M$ back to the hyperbolic metric, up to rescaling, after finitely many surgeries.  An important point to check is that the minimal surfaces we consider avoid the surgery regions. The regions that are removed in surgery resemble long thin tusks capped at the end by a (potentially much larger) region homeomorphic to a ball.  A simple argument using the monotonicity formula guarantees that the minimal surfaces we consider do not enter these regions.  To finish the proof we use the fact that normalized Ricci flow takes small perturbations of the constant curvature metric back to the constant curvature metric at an exponential rate of convergence. The exponent for the convergence is determined by the linearization of the RHS of the normalized Ricci flow we consider.  We use the work of Knopf-Young \cite{knopfccf} and the assumption that $M$ admits no infinitesimal conformally flat deformations to bound the spectrum of the linearization and control the rate of convergence to the hyperbolic metric.

\subsection{} The strategy of obtaining bounds on areas of surfaces via Ricci flow is by no means original to this paper (\cite{bben}, \cite{coldminiextinction},\cite{marquesnevesduke}.)  Conceivably other evolution equations for metrics could be used to similar effect, although we only know of the paper by Ambrozio-Montezuma \cite{ambmont}, which studied the Simon-Smith widths of conformal classes of metrics on $S^3$ via Yamabe flow.   

The work of Bray-Brendle-Eichmair-Neves \cite{bben}  is the closest thing to a positive curvature analogue of the results here.  For closed 3-manifolds with a positive lower scalar curvature bound, they prove sharp upper bounds for the area of an embedded projective plane.  They also show that the case of equality implies that the 3-manifold in question had to be isometric to the standard $\mathbb{RP}^3$.  Unlike this paper, however, they only use Ricci flow in characterizing the case of equality, where moreover they only use the short-time existence.

The paper by Marques-Neves  \cite{marquesnevesduke} on the other hand, which obtained sharp area bounds for min-max minimal surfaces in 3-manifolds with non-negative Ricci curvature under positive scalar curvature lower bounds, did use the long-time behavior of Ricci flow. (But see \cite{songembness}, which removes the non-negative Ricci curvature assumption and  uses only the short-time existence.) Both \cite{bben} and \cite{marquesnevesduke} were important sources of inspiration for this paper.

\subsection{} \label{relatedresults} We now describe some other results related to this paper.  In \cite{andersoncanmets}, Anderson gives an argument to show how Perelman's work allows one to compute the Yamabe invariant of a hyperbolic 3-manifold $M$.  As in this paper, he starts with a general metric and evolves it under Ricci flow with surgery while keeping track of how scalar curvature and the relevant geometric quantities are changing. A simple consequence of the computation of the Yamabe invariant is that any metric on $M$ that is not hyperbolic and that has scalar curvature greater than or equal to that of the hyperbolic metric must have volume greater than that of the hyperbolic metric.  (This is actually conjectured to hold in all dimensions; see \cite{guthlargeballs} for further discussion and some results in that direction.)  

A landmark comparison theorem for negative scalar curvature lower bounds  is the hyperbolic positive mass theorem, proved by Min-Oo in the spin case and Andersson-Cai-Galloway in dimension less than 8.  It implies as a special case that a metric on $\mathbb{H}^n$ with scalar curvature greater than or equal to $-n(n-1)$ that is equal to the standard hyperbolic metric outside of a compact set must actually be isometric to the standard $\mathbb{H}^n$. This special case also follows from recent work of Li \cite{chaolipolyhyp}, who proved a comparison theorem for polyhedra in manifolds with a negative lower scalar curvature bound. 

Finally we mention the paper by Ache-Viaclovsky \cite{achvia15}, which also considers hyperbolic 3-manifolds that admit no infinitesimal conformally flat deformations in a geometric analysis context.

\subsection{Acknowledgments} 

I thank Chao Li for a helpful conversation and for telling me about the Bamler paper \cite{bamlermain}. I thank Antoine Song for a helpful conversation. I also thank my advisor Fernando Coda Marques for a helpful conversation related to this paper and for his support.  I am indebted to Andre Neves for bringing to my attention an error in the first version of this paper and for his efforts to help me fix it.

	\section{Proof of Theorem \ref{varcurvintro}} \label{first} 
Let $M$ be a closed hyperbolic 3-manifold and let $\Sigma$ be a closed immersed $\pi_1$-injective surface in $M$. Let $\mathcal{F}_\Sigma$ denote the set of all immersed surfaces in $M$ homotopic to the inclusion of $\Sigma$ in $M$.  For a Riemannian metric $g$ on $M$, define
\[
\mathcal{A}_{\Sigma}(M,g) = \inf \{ \text{Area}_g(\Sigma) : \Sigma \in \mathcal{F}_\Sigma\}.  
\]

   
In this section we will first prove Theorem \ref{varcurvintro} in the special case that $\Sigma$ is totally geodesic in $(M,g_{hyp})$. For $(M,g_{hyp})$ to satisfy the hypothesis of infinitesimal conformal rigidity, $\Sigma$ must necessarily fail to be embedded \cite{JM87}.  The author does not know whether there are infinitesimally conformally rigid closed hyperbolic 3-manifolds that contain immersed totally geodesic surfaces, so this case might well be vacuous.  We give the proof of this case, however, since it is illustrative and contains all of the ideas for the proof of the general case of $\Sigma$ almost totally geodesic. In this special case we are also able to prove that $g$ is isometric to $g_{hyp}$ in the event that equality in the inequality is realized and
\[
   \mathcal{A}_{\Sigma}(M,g)= -2\pi \chi(\Sigma)  = \mathcal{A}_{\Sigma}(M,g_{hyp}), 
   \]
which anticipates the rigidity in the case of equality in Theorem \ref{cmnscalintro}.

First, we claim that we can assume that $\Sigma$ is two-sided (has trivial normal bundle in $M$.)  If $\Sigma$ is not two-sided, then we can find an immersion of a double cover $\Sigma'$ of $\Sigma$ that is.  The claim then follows from the fact that homotopies of $\Sigma$ lift to homotopies of $\Sigma'$ and the area of an immersion of $\Sigma$ is half that of the immersion of $\Sigma'$ that lifts it.  From now on we assume that $\Sigma$ is two-sided.    
	
The proof  will be based on the Ricci flow with surgery. By work of Hamilton-Perelman, we know that Ricci flow with surgery starting at $(M,g)$ converges to the hyperbolic metric after rescaling the metrics to have volume that of the hyperbolic metric.  To prove the inequality, we study how $\mathcal{A}_\Sigma(M,g)$ evolves under Ricci flow in comparison to $\mathcal{A}_\Sigma(M,g_{hyp})$. 

\begin{proof} 
	There exists a real number $T$ and a family of metrics $g(t)$ such that $g(t)$ satisfies the Ricci flow equation 
	\begin{equation} \label{ricflow} 
	\frac{\partial}{\partial t} g(t) = -2 \text{Ric}_{g(t)}   
	\end{equation} 
	\noindent for $t \in [0,T)$ (\cite{deturcktrick},\cite{hamshorttime}.)  It will be important for us that the scalar curvature satisfies the evolution equation  
	\begin{equation} \label{scalevol} 
		\frac{\partial}{\partial t} R_{g(t)} = \Delta R_{g(t)} + 2 |\text{Ric}_{g(t)}|^2. 
	\end{equation}  
	
	If we set $g=g_{hyp}$, then $g_{hyp}(t) = (1+4t)g_{hyp}$ solves (\ref{ricflow}) with $T=\infty$.  If $T< \infty$, then high curvature regions develop as $t \nearrow T$.  We can perform surgery to cut out the high curvature regions in a standard and controlled way, so that bounds on all relevant geometric quantities continue to hold after the surgery, and we can then 	restart Ricci flow and repeat. The work of Perelman implies that this can be done in such a way that only finitely many surgeries occur total, which we explain in more detail below.  In the case of our $(M,g)$, which cannot be written as a non-trivial connect sum, every surgery removes an inessential neck bounding a 3-ball.

	\subsection{$\Sigma$  Totally Geodesic Case 1: No Surgeries} 
	
	To start out, we're going to give the proof for Riemannian metrics $g$ such that there exist $\{g(t): t \in [0,\infty)\}$ that satisfy (\ref{ricflow}) for all time and converge to the hyperbolic metric on $M$ after rescaling the metrics to have volume that of the hyperbolic metric. 
	
	Let $\Sigma$ be a closed immersed totally geodesic surface in $(M,g_{hyp})$. We take $\Sigma_t$ to be a closed immersed minimal surface in $(M,g(t))$ that minimizes $g(t)$-area over all surfaces homotopic to the inclusion of $\Sigma$ in $M$. There exists at least one such surface by \cite{sacksuhl} or \cite{schoenyauincompressible}. For fixed $t_0$ and $\text{Area}_{\Sigma_{t}}(s)$ the area of $\Sigma_t$ in the metric $g(s)$, we have that
	\begin{equation} \label{trric}  
	\frac{d}{ds}(\text{Area}_{\Sigma_{t_0}}(s))|_{s=t_0} = - \int_{\Sigma_{t_0}} \text{Ric}_{g(t_0)}(E_1,E_1) +  \text{Ric}_{g(t_0)}(E_2,E_2) d \mu_{g(t_0)},
	\end{equation}    	 
	where $(E_1,E_2)$ is a local orthonormal frame on $\Sigma_t$ for $g(t_0)$ and $d \mu_{g(t_0)}$ is the area form for the metric on $\Sigma_{t_0}$ induced by $g(t_0)$.  Recall that we assume by the comment following Theorem \ref{varcurvintro}  that $\Sigma$ is two-sided. If $\nu$ is the unit normal vector to $\Sigma$ we can write the integrand of the above equation as
	
	\vspace{-6mm} 
	\begin{align*} 
	 &R(E_1,\nu,E_1,\nu) + R(E_2,\nu,E_2,\nu) + 2R(E_1,E_2,E_1,E_2)  \\ &=  \frac{1}{2} R + R(E_1,E_2,E_1,E_2)    \\
	&=   \frac{1}{2} R + K + \frac{1}{2} |A|^2,   
	\end{align*}  	
	\noindent where $R(\cdot,\cdot,\cdot,\cdot)$ and $R$ are respectively the curvature tensor and the scalar curvature of $(M,g(t))$, $K$ is the Gauss curvature of $\Sigma_{t_0}$, and $A$ is the second fundamental form of $\Sigma$.  To get the last equality, we're using the Gauss equation and the fact that $\Sigma_{t_0}$ is minimal. Using Gauss-Bonnet and the fact that $|A|^2 \geq 0$, we then have that 
	\begin{equation} \label{darea} 
	 	\frac{d}{ds}(\text{Area}_{\Sigma_{t_0}}(s))|_{s=t_0} \leq \int_{\Sigma_{t_0}} -\frac{1}{2} R   d \mu_{g(t_0)} - 2\pi \chi( \Sigma).  
	\end{equation} 
	
	Let $R_{min}(t)$ be the  minimum of the scalar curvature of $g(t)$.  Then since $R_{min}(0) = -6$, we have by comparing with Ricci flow starting at $g_{hyp}$ and the parabolic maximum principle that 
	\begin{equation} \label{minest} 
	R_{min}(t) \geq \frac{-6}{(4t+1)}.  
	\end{equation} 

	

	 The inequalities (\ref{darea}) and (\ref{minest}) imply that 
	\begin{equation} \label{spec} 
	\frac{d}{ds}(\text{Area}_{\Sigma_{t_0}}(s))|_{s=t_0} \leq \frac{3 \text{Area}_{\Sigma_{t_0}}(t_0)}{4t_0+1} -  2\pi \chi( \Sigma).  
	\end{equation}   
	
	Set $\mathcal{A}(t)= \mathcal{A}_{\Sigma}(M,g(t))$.  If $g(0)=g_{hyp}$, we set  $\mathcal{A}_{hyp}(t)= \mathcal{A}_{\Sigma}(M,g(t))$.  By the same argument as the proof of Lemma 9 of \cite{bben}, the function $\mathcal{A}(t)$ is Lipschitz and thus differentiable almost everywhere.  	For any fixed $t_0$, 
	\[
	\mathcal{A}(t)\leq \text{Area}_{\Sigma_{t_0}}(t). 
	\]

\noindent We also have that 	
	\[
	\mathcal{A}(t_0) =  \text{Area}_{\Sigma_{t_0}}(t_0). 
	\]
	
\noindent 	If $\mathcal{A}(t)$ is differentiable at $t_0$, we therefore have 
	\[
		\mathcal{A}'(t_0)= \text{Area}_{\Sigma_{t_0}}'(t_0).  
	\] 
	
	\noindent For almost every $t$, $\mathcal{A}(t)$ thus by (\ref{spec}) satisfies the inequality 
	\begin{equation} \label{keyineq} 
	\mathcal{A}'(t) \leq \frac{3\mathcal{A}(t)}{4t+1} -  2\pi \chi( \Sigma).	
	\end{equation} 
	
\noindent	If the metric we started with was hyperbolic, then we would have that 
	\begin{equation} \label{ineqmodel} 
			\mathcal{A}_{hyp}'(t) = \frac{3\mathcal{A}_{hyp}(t)}{4t+1} -  2\pi \chi( \Sigma). 
	\end{equation} 
	 Therefore, since $\mathcal{A}(0)$ is less than or equal to $\mathcal{A}_{hyp}(0)$, it follows by taking the difference of (\ref{ineqmodel}) and (\ref{keyineq}) and solving explicitly the ODE 
	 \[
	 y'(t) = \frac{3y(t)}{4t+1} 
	 \]	 
	 obtained by replacing the inequality sign with an equality sign that
	 \begin{equation} \label{difference}
	 \mathcal{A}_{hyp}(t) - \mathcal{A}(t) \geq (4t + 1)^{\frac{3}{4}} (\mathcal{A}_{hyp}(0) - \mathcal{A}(0)). 
	 \end{equation}
	To get (\ref{difference}), we are using the fact that $\mathcal{A}$ is, as a Lipschitz function, absolutely continuous, and so the fundamental theorem of calculus can be applied.
	 


	 As $t$ tends to infinity, the metrics $g(t)$ converge to the hyperbolic metric when rescaled to have volume that of the hyperbolic metric.  Denote these rescaled metrics by $\hat{g}(t)$.  Then $\hat{g}(t) - g_{hyp}$ can be made pointwise arbitrarily small in any $C^k$ norm by taking $t$ large enough. 

	 Knopf-Young \cite{knopfccf} consider a normalized Ricci flow, which they call KNRF, defined by the following equation  
	 \begin{equation} \label{knrf} 
	 \frac{\partial}{\partial t} g(t) = -2 \text{Ric}_{g(t)} - 4g(t). 
	 \end{equation} 
	 
	Note that the hyperbolic metric is fixed by KNRF. Solutions to KNRF differ from solutions to usual Ricci flow by rescalings in space and time. If $\tilde{g}(\tilde{t})$ is a solution to (\ref{knrf}), then the time variable $t$ for the corresponding Ricci flow $g(t)$ is equal to $\frac{1}{4} \log (1+4\tilde{t})$. For metrics $\overline{g}$ in a sufficiently small neighborhood $U$ of the hyperbolic metric in the $C^{2,\alpha}$-norm, Knopf-Young show that the solution to (\ref{knrf}) starting at $\overline{g}$ converges back to the hyperbolic metric exponential fast. They do so by bounding the spectrum of the linearization $A_{g_{hyp}}$ at the hyperbolic metric of the RHS of (\ref{knrf}), which they compute to be 
	 \begin{equation} 
	 A_{g_{hyp}}(h) = \Delta h - 2H g_{hyp} + 2h,
	 \end{equation} 
	 \noindent where $h$ is a symmetric $(2,0)$ tensor and $H=g^{ij}h_{ij}$ is the trace of $h$.  The linearized operator $A_{g_{hyp}}$ is self-adjoint and elliptic with discrete spectrum contained in an interval $(-\infty, \omega]$. If $\omega<0$ and $\overline{g}(t)$ is the KNRF starting at $\overline{g}$ contained in the $C^{2,\alpha}$ neighborhood described above, then the $C^0$ norm $\overline{g}(t)-g_{hyp}$ is bounded above by a constant times $e^{\omega_0t}$ for any $\omega_0>\omega$.  Using a Bochner formula for symmetric $(2,0)$ tensors due to Koiso, Knopf-Young show that $\omega$ can be taken to be equal to $-1$. 
	 
	  Let $h$ be an eigentensor with eigenvalue $-1$. Then it follows from \cite{knopfccf}[Section 5] that $h$ must be trace-free and Codazzi, where recall that a symmetric $(2,0)$ tensor $T$ is Codazzi if its covariant derivative tensor is also symmetric-- i.e.,  
	  \[
	  \nabla_X T(Y,Z) = \nabla_{Y}T(X,Z).  
	  \]  
	 
	 Trace-free Codazzi tensors over a hyperbolic 3-manifold correspond to infinitesimal conformally flat deformations of the hyperbolic metric (see, e.g., the introduction of \cite{beig97}.)  The space of infinitesimal conformally flat deformations is given by the cohomology group 
	 \begin{equation} \label{tang} 
	 H^1(\pi_1(M),\text{Ad}), 
	 \end{equation} 
	 \noindent where $\text{Ad}$ is the Adjoint representation of $\pi_1(M)\subset SO(3,1)$ on $so(4,1)$ via the inclusion $so(3,1) \hookrightarrow so(4,1)$.  As mentioned in the introduction, Kapovich \cite{kap94}  gave examples of closed hyperbolic 3-manifolds for which this cohomology group vanishes. For such infinitesimally rigid $M$ there can be no eigentensor $h$ as above with eigenvalue $-1$, and so $\omega$ can be taken to be less than $-1$.  On the other hand, we point out that if (\ref{tang}) is nonzero, which happens for example whenever $(M,g_{hyp})$ contains an embedded totally geodesic surface \cite{JM87}, then there will be eigentensors with eigenvalue $-1$, and so the argument below would fail.

	 Now assume that $(M,g_{hyp})$ admits no infinitesimal conformal deformations.  Let $\hat{g}(t)$ be the rescaling of our initial Ricci flow $g(t)$ to have volume equal to that of $g_{hyp}$ at all times.  Let $\tilde{g}(\tilde{t})$ be the KNRF with $\tilde{g}(0)=g(0)$. Since $\tilde{g}(\tilde{t})$ $C^0$-converges to $g_{hyp}$ at a rate of $e^{\omega \tilde{t}}$ for some $\omega<-1$, we have that $\hat{g}(t)$ $C^0$-converges to the hyperbolic metric at a rate of $t^{\omega/4}$, since $\tilde{t}= \frac{1}{4}\log(1+4t)$.  This is because after some finite time $\hat{g}(t)$ will be contained in the $C^{2,\alpha}$ neighborhood $U$ of $g_{hyp}$ described above.

	Dividing Equation (\ref{difference}) by $(Vol(g(t))/Vol(g_{hyp}))^{2/3}$, which is less than or equal to a constant times $t$, 
	we have that for some positive constant $C_1$  
	\[
	 \mathcal{A}_{\Sigma}(M,g_{hyp})-\mathcal{A}_{\Sigma}(M,\hat{g}(t)) \geq \frac{C_1}{t^{1/4}}. 
	\]
    To obtain this inequality we have also used the fact that $Vol(g(t))\geq Vol(g_{hyp}(t))$, which follows from the theorem in \cite{andersoncanmets} described in subsection \ref{relatedresults} of the introduction, since the scalar curvature of $g(t)$ is greater than that of $g_{hyp}(t)$ at all times $t$.

	  But taking $t$ large  and using $\hat{g}(t)$'s convergence to the hyperbolic metric at a rate of $t^{\omega/4}$ gives a contradiction, since this implies that for some constant $C_2$ 
	\begin{equation} \label{endcontr} 
	|\mathcal{A}_{\Sigma}(M,g_{hyp})-\mathcal{A}_{\Sigma}(M,\hat{g}(t))| \leq C_2t^{\omega/4},  
	\end{equation} 
	\noindent and $\omega<-1$.

	\vspace{.5mm} 
	\noindent \subsubsection{Case of Equality} \label{caseofequality}  Now suppose that for some metric $g$ with scalar curvature greater than or equal to $-6$, 
	\[
	 \mathcal{A}_{\Sigma}(M,g)=\mathcal{A}_{\Sigma}(M,g_{hyp}).
	 \]    
	 Then we will show, following the argument for the equality case of the main result of \cite{bben}, that $g$ is isometric to $g_{hyp}$.  First, we run the Ricci flow for a small interval of time $[0,\epsilon)$. We write the evolution equation for scalar curvature under Ricci flow as
	\begin{equation} \label{modevol} 
		\frac{\partial}{\partial t} R_{g(t)} =\Delta R_{g(t)} + \frac{2}{3} R_{g(t)}^2 + 2 |\mathring{Ric}_{g(t)}|^2.  
	\end{equation} 
	


 By (\ref{difference}), $\mathcal{A}_{\Sigma}(M,g(t))\leq \mathcal{A}_{\Sigma}(M,g_{hyp}(t))$ for $t\in[0,\epsilon)$, so by 
 what we have shown above we must have that the minimum of $R_{g(t)}$ is identically equal to $\frac{-6}{1+4t}$ on $[0,\epsilon)$. Consequently, by the strong maximum principle, $R_{g(t)}$ is identically equal to $\frac{-6}{1+4t}$ which implies that $\mathring{Ric}_{g(t)}$ is identically zero on $[0,\epsilon)$. The metric $g$ must then have been Einstein and, since we are in three dimensions, have had constant sectional curvature $-1$.  Since by Mostow rigidity the constant curvature $-1$ metric on $M$ is unique up to isometry, this implies that $g$ is isometric to $g_{hyp}$.

	\end{proof} 
	
	\subsection{$\Sigma$ Totally Geodesic Case 2: Finitely Many Surgeries}
	
	We now modify the proof of the previous case in the event that Ricci flow $(M,g(t))$ starting from the initial metric  develops singularities. First we describe the general strategy. At a singularity, long tusks are forming that up to rescaling are nearly isometric to $(S^2, g_{round}) \times (-\frac{1}{\epsilon},\frac{1}{\epsilon})$ along most of their length, for $\epsilon$ small.  In performing surgery, we chop off part of the tusk.  In order to implement the approach of the previous case, we need to rule out the possibility that an area-minimizing surface in the homotopy class of $\Sigma$ enters the part of the tusk we chop off in surgery, which we accomplish by using the monotonicity formula.  This proves that $\mathcal{A}_{\Sigma}(M,g(t))$ is Lipschitz continuous at the surgery times.  Bounds on all relevant geometric quantities--- most importantly for us the pointwise lower bound on the scalar curvature--- continue to hold after surgery, so from here the argument can proceed as in the first case.  
	
	We now explain Bamler's set-up, which we copy from his paper \cite{bamlermain} nearly verbatim. We don't actually use any of his results from that paper, only his definitions and formulations of Perelman's results, which are concise and convenient for our purposes. To start we give the definition of Ricci flow with surgery that we will use.     
	
	\begin{defn} [Ricci flow with surgery] 
		Consider a time interval $I\subset \mathbb{R}$.  Let $T^1<T^2<...$ be times in the interior of $I$ which form a possibly infinite  but discrete subset of $\mathbb{R}$ and divide $I$ into the intervals 
		\[
		I^1 = I \cap (-\infty, T^1), I^2 = [T^1,T^2), ...
		\]
		and $I^{k+1}= I \cap [T^k, \infty)$ if there are only finitely many $T^i$'s.  Consider Ricci flows $(M^1 \times I^1,g_t^1)$, $(M^2 \times I^2, g_t^2)$,... on 3-manifolds $M^1$, $M^2$,... 
		Assume that the metrics $g_t^i$ converge smoothly as $t \nearrow T^i$ to a Riemannian metric $g_{T^i}$ on $M$, and let 
		\[
		U_{-}^i \subset M^i \hspace{2mm} \text{and} \hspace{2mm} U_{+}^i \subset M^{i+1} 
		\]
		
		\noindent be open subsets such that there are isometries 
		\[
		\Phi^i:(U_{-}^i, g_{T^i}^i) \rightarrow (U_{+}^i,g_{T^i}^{i+1}), \hspace{4mm} (\Phi^i)^* g_{T^i}^{i+1}|_{U_{+}^i} = g_{T^i}^i|_{U_{-}^i}.
		\]
		
		We assume that we never have $U_{-}^i = M^i$ and $U_{+}^i= M^{i+1}$. 
	  We also assume that every component of $M^{i+1}$ contains a point of $U_+^i$.  We call $\mathcal{M} = ((T^i)_i, (M^i \times I^i, g_t^i)_i,(U_{\pm}^i)_i, (\Phi^i)_i)$ a \textit{Ricci flow with surgery} on the time interval $I$, and we call $T^1$,$T^2$,... surgery times.  
		
		If $t\in I^i$, then $(\mathcal{M}(t),g(t))= (M^i \times \{t\}, g_t^i)$ is called the time-$t$ slice of $\mathcal{M}$.  For $t=T^i$, we define the (presurgery) time $T^{i-}$-slice to be $(\mathcal{M}(T^{i-}), g(T^{i-}))= (M^i \times  \{T^i\},g_{T^i}^i)$.  The points $M^i \times \{T^i\} \setminus  U_-^i \times \{T^i\}$ are called presurgery points and the points $M^{i+1} \times \{T^i\}\setminus U_+^i \times \{T^i\}$ are called surgery points.  We will call a point that is not a presurgery point a non-presurgery point.   
		\end{defn} 
	\begin{rem} 
		This definition of Ricci flow with surgery is a slight specialization of the one from \cite{bamlermain}, since the only surgeries that occur for our $M$ pinch off a topologically trivial capped neck.

		\end{rem} 
	
		\begin{defn}[$\epsilon$-neck] 
		Let $\epsilon>0$, and consider $U \subset (M,g)$.  Then we say $U$ is an $\epsilon$\textit{-neck} if there is a diffeomorphism $\Phi: S^2 \times (-\frac{1}{\epsilon},\frac{1}{\epsilon}) \rightarrow U$  and a $\lambda>0$ such that $|\lambda^{-2}  \Phi^* g - g_{S^2 \times (-\frac{1}{\epsilon},\frac{1}{\epsilon})}|_{C^{\lceil \epsilon^{-1} \rceil}}<\epsilon $, where $g_{S^2 \times (-\frac{1}{\epsilon},\frac{1}{\epsilon})}$ is the standard metric on $S^2 \times (-\frac{1}{\epsilon},\frac{1}{\epsilon})$ with constant scalar curvature 2.  We say that $x$ is a \textit{center} of $U$ if $x \in \Phi(S^2 \times {0})$. 			
	\end{defn} 
	
It follows from the results of Perelman stated in \cite{bamlermain} that for any metric $g$ on $M$ and any $\epsilon>0$ to be specified later that we choose at the start, there is a Ricci flow with surgery $g(t)$ starting at $g$ and defined for all time that satisfies the following.  At every surgery time $T^i$, the components of $\mathcal{M}(T^{i-}) \setminus U_{-}^i$ are homeomorphic to $\mathbb{D}^3$, and the subset $\mathcal{M}(T^i) \setminus   U_+^i$ is a disjoint union $D_1^i \cup .. \cup D_{m_i}^i$ of homeomorphic copies of $\mathbb{D}^3$.  Moreover, for every $D_j^i$, we can take the points on the boundary of $U_-^i$ in $\mathcal{M}(T^{i-})$ corresponding to $\partial D_j^i$ to be centers of $\epsilon$-necks for the $\epsilon$ chosen at the start.  This follows from Proposition 3.4 of \cite{bamlermain}, the definition of Ricci flow with surgery with $\delta(t)$-precise cutoff, and the fact that we can take $\delta(t)$ to be pointwise smaller than any positive number (in this case, $\epsilon$.) 

Moreover, for our specific $M$ only finitely many surgeries occur.  This is because the $g(t)$, normalized to have volume that of the hyperbolic metric, smoothly converge to the hyperbolic metric as $t$ tends to infinity. This is explained in \cite{andersoncanmets}[pg. 132] and \cite{calegaririccinotes}[pg. 60]. Since only finitely many surgeries occur on any finite time interval, and no surgeries occur when the normalized metric is sufficiently close to the hyperbolic metric, only finitely many surgeries occur total.


	
	Let $\Sigma_{T^{i-}}$ be an immersed minimal surface in $\mathcal{M}(T^{i-})= (M,g_{T^{i-}})$ that realizes the minimum $g_{T^{i-}}$-area over all immersed surfaces homotopic to $\Sigma$. We claim that $\Sigma_{T^{i-}}$ does not intersect $\mathcal{M}(T^{i-}) \setminus U_{-}^i$, which is diffeomorphic to a union of balls, provided that the $\epsilon$ chosen at the start was taken sufficiently small.  Any point $p$ on the boundary of this region is the center of an $\epsilon$-neck $N$. 
	\[
	\Phi: S^2 \times (-\frac{1}{\epsilon}, \frac{1}{\epsilon}) \rightarrow N 
	\]
	such that $p\in \Phi(S^2 \times \{0\})$ and  for some $\lambda>0$
	\begin{equation} \label{closeness} 
	|\lambda^{-2}  \Phi^* g_{T^{i-}} - g_{S^2 \times (-\frac{1}{\epsilon},\frac{1}{\epsilon})}|_{C^{\lceil \epsilon^{-1} \rceil}}<\epsilon. 
	\end{equation} 
	
	 
	  Now assume for contradiction that $\Sigma_{T^{i-}}$ passes through $p$.  It cannot be the case that $\Sigma_{T^{i-}}$ is contained in $N \cup  \mathcal{M}(T^{i-}) \setminus U_{-}^i$, so $\Sigma_{T^{i-}}$ must intersect $\Phi(S^2 \times \{x \})$, for $x$ slightly greater than $-\frac{1}{\epsilon}$.  We can perturb $\Phi(S^2 \times \{x \})$ slightly to an embedded sphere $S$ that intersects $\Sigma_{T^{i-}}$ transversely in a union of circles. These circles necessarily bound disks in $\Sigma_{T^{i-}}$, so let $D_p$ be a disk or an annulus in $\Sigma_{T^{i-}}$ that passes through $p$.  
	 
	 On the one hand, $\partial{D_p}$ can be filled in by a region in $S$ with area at most roughly $4\pi\lambda^2$.  On the other hand, $D_p$ must intersect every cross-section $\Phi(S^2 \times \{y \})$ for $-\frac{1}{\epsilon}<y<0$.  By the monotonicity formula, there is a universal constant $c$ such that for $y \in (-\frac{1}{\epsilon}+1/2, -1/2)$, 	 
	 \[
	 \text{Area}(D_p \cap \Phi( S^2 \times (y-1/2,y+1/2)) > c    \lambda^{2}.    
	 \]
	 It follows by choosing disjoint unit intervals in $(-\frac{1}{\epsilon}+1/2, -1/2)$ that if we chose $\epsilon$ such that $\frac{1}{\epsilon}$ is greater than $5\pi/c$, then the area of $D_p$ will be larger than $(4\pi+1)\lambda^2$.  By cutting out $D_p$ and gluing in a region of  $S$ we could then produce a surface homotopic to $\Sigma_{T^{i-}}$ but with smaller area, which is a contradiction.  
	 
	 Since $U_{-}^i$ is unchanged by surgery and we have an isometry $U_{-}^i \rightarrow U_{+}^i$, the same argument shows that a $\Sigma_{T^{i+}}$ that minimizes area over surfaces homotopic to $\Sigma$ in the surgered manifold $\mathcal{M}^{T^{i+}}$ is contained in  $U_{+}^i$.  This proves that the function $\mathcal{A}(t)$,  defined as the minimum area of a surface homotopic to $\Sigma$ in $\mathcal{M}(t)$, is well-defined.  Since it is also Lipschitz, and since the estimate (\ref{minest})  for the minimum of scalar curvature also holds for Ricci flow with surgery, the arguments of the previous section apply. This proves Theorem \ref{varcurvintro} in the case that $\Sigma$ is totally geodesic in the hyperbolic metric, in the event that singularities occur in Ricci flow starting at the initial metric.   
	 
	\noindent \subsection{Proof for General $\Sigma$ }  
	
 Fix a metric $g$ with scalar curvature greater than or equal to $-6$, and let $g(t)$ be a Ricci flow with surgery starting from $g$.  The inequality in Theorem \ref{varcurvintro} can be rearranged to 
	\begin{equation} \label{startineq} 
	\mathcal{A}_{\Sigma}(M,g_{hyp}) - \mathcal{A}_{\Sigma}(M,g) \leq \delta \mathcal{A}_{\Sigma}(M,g_{hyp}).
	\end{equation}

	\noindent Fixing a surface $\Sigma$ as in the theorem and with notation as above, we assume that for some $\delta>0$ 
	\[
	 \mathcal{A}_{hyp}(0) - \mathcal{A}(0) > \delta \mathcal{A}_{hyp}(0).  
	\]

	\noindent Since $\Sigma$ is not totally geodesic, in place of (\ref{ineqmodel}) we have the equality	
	\begin{equation} \label{ineqmodel2} 
	\mathcal{A}_{hyp}'(t) = (1- \epsilon_{\Sigma}) \left( \frac{3\mathcal{A}_{hyp}(t)}{4t+1}\right)  -  2\pi \chi( \Sigma). 
	\end{equation} 
	
	\noindent where recall that $\epsilon_{\Sigma}$ is the average of $|A|^2$ over $\Sigma$.
Taking the difference with (\ref{keyineq}), we obtain 
\begin{equation} \label{differror} 
(\mathcal{A}_{hyp}(t)- \mathcal{A}(t))' \geq  \left( \frac{3(\mathcal{A}_{hyp}(t)- \mathcal{A}(t))}{4t+1}\right) - \frac{3\epsilon_{\Sigma}}{4t+1} \left(\mathcal{A}_{hyp}(t)\right)  .  
\end{equation} 

\noindent Since $\mathcal{A}_{hyp}(t)= (1+4t) \mathcal{A}_{hyp}(0)$, integrating gives 
\begin{align} 
\mathcal{A}_{hyp}(t)- \mathcal{A}(t) &\geq (4t + 1)^{\frac{3}{4}} (\mathcal{A}_{hyp}(0) - \mathcal{A}(0)) - 3 \epsilon_{\Sigma}  \mathcal{A}_{hyp}(0)t   \\
\label{finalcontr} &> (4t + 1)^{\frac{3}{4}} (\delta \mathcal{A}_{hyp}(0)) - 3 \epsilon_{\Sigma}  \mathcal{A}_{hyp}(0)t.
\end{align} 

Dividing by $(Vol(g(t))/Vol(g_{hyp}))^{2/3}$ as above, which grows at most  linearly, and using the fact that $\hat{g}(t)$ converges to $g_{hyp}$ at a rate of at least $t^\frac{\omega}{4}$ for $\omega<-1$ as $t$ tends to infinity, which implies the inequality (\ref{endcontr}), we obtain a contradiction if $\epsilon_{\Sigma}$  is sufficiently small.  This shows that for every $\delta>0$, (\ref{startineq}) holds for $\epsilon_{\Sigma}$ sufficiently small.  

	\section{Counting with a Lower Scalar Curvature Bound} \label{second} 
	
In this section, we prove Theorem \ref{cmnscalintro}. Assume that $(M,g)$ is as in Theorem \ref{cmnscalintro}, with scalar curvature greater than or equal to $-6$.  By \cite{sacksuhl} or \cite{schoenyauincompressible}, every element of $S_{\epsilon}(M)$ is realized by an immersed $g_{hyp}$-minimal surface which has principal curvatures at most $C\log(1+\epsilon)$, for some universal constant $C$ \cite{seppi}.  For every $\delta>0$, Theorem \ref{varcurvintro} implies that for $\epsilon$ sufficiently small 
	\begin{multline} 
	\# \{ \text{Area}_g (S) \leq 4 \pi(L-1) : S \in S_{\epsilon}(M) \} \\ 	
	\leq \# \{\text{Area}_{g_{hyp}} (S) \leq 	4 \pi(L-1)/(1-\delta) : S \in S_{\epsilon}(M)\}.   
	\end{multline}
	\noindent Dividing by $L \log L$, sending $L$ to infinity, and sending $\delta$ to $0$ proves that $E(g) \leq E(g_{hyp})$.
	
	Now assume $E(g) = E(g_{hyp})$.  We follow the same strategy as in \ref{caseofequality}. Let $g(t)$ be a Ricci flow starting at $g$ for $t$ contained in a small interval of time $[0,\epsilon')$. We claim  that $E(g(t)) \geq E(g_{hyp}(t))$, where $g_{hyp}(t)$ is the Ricci flow starting at $g_{hyp}$.  To see this, we claim that for every $\delta>0$ 
	\begin{equation} \label{smallpower} 
	\text{Area}_{g(t)}(S) \leq (1+4t)^{1+\delta} \text{Area}_{g}(S),
	\end{equation} 
	\noindent provided $\epsilon$ was taken sufficiently small and $S\in S_{\epsilon}(M)$.  We have that for any $\delta'>0$, 
	\[
	\text{Area}_g(S) \geq (1-\delta') (-2 \pi \chi(S))
	\]
provided $\epsilon$ was taken small enough. Hence for $\epsilon$ small enough and $\delta'$ taken small enough relative to $\delta$, we have that $(1+4t)^{1+\delta} \text{Area}_{g}(S)$ is a super-solution of the ODE 
	\begin{equation} \label{inequallity} 
	\mathcal{A}'(t) = \frac{3\mathcal{A}(t)}{4t+1} -  2\pi \chi( \Sigma).  
	\end{equation} 
	If $S$ were totally geodesic and had area equal to $-2 \pi \chi(S)$, $(1+4t) \text{Area}_{g}(S)$ would satisfy this ODE.  Inequality (\ref{keyineq}) says exactly that  $\text{Area}_{g(t)}(S)$ is a sub-solution of (\ref{inequallity}.) The inequality (\ref{smallpower}) now follows.

	We also have that
	\[
	\text{Area}_{g_{hyp}(t)}(S) = (1+4t) \text{Area}_{g_{hyp}}(S).  
	\]
	\noindent Therefore, since $E(g(0))=E(g_{hyp}(0))$, sending $\delta$ to $0$ we conclude that $E(g(t)) \geq E(g_{hyp}(t))$, so by what we showed at the start of this section, this implies that the minimum of the scalar curvature of $g(t)$ is less than or equal to the scalar curvature of $g_{hyp}(t)$. From here we can reason as we did in \ref{caseofequality} to conclude that $g$ is isometric to $g_{hyp}$.

		\bibliography{bibliography}{}

\begin{thebibliography}{BBEN10}

\bibitem[AM18]{ambmont}
Lucas Ambrozio and Rafael Montezuma.
\newblock On the width of unit volume three-spheres.
\newblock {\em arXiv:1809.03638 [math]}, 2018.

\bibitem[And06]{andersoncanmets}
Michael~T. Anderson.
\newblock Canonical metrics on 3-manifolds and 4-manifolds.
\newblock {\em Asian J. Math.}, 10(1):127--163, 2006.

\bibitem[AV15]{achvia15}
Antonio~G. Ache and Jeff~A. Viaclovsky.
\newblock Asymptotics of the self-dual deformation complex.
\newblock {\em J. Geom. Anal.}, 25(2):951--1000, 2015.

\bibitem[Bam17]{bamlermain}
Richard~H. Bamler.
\newblock The long-time behavior of 3-dimensional {R}icci flow on certain
  topologies.
\newblock {\em J. Reine Angew. Math.}, 725:183--215, 2017.

\bibitem[BBEN10]{bben}
H.~Bray, S.~Brendle, M.~Eichmair, and A.~Neves.
\newblock Area-minimizing projective planes in 3-manifolds.
\newblock {\em Comm. Pure Appl. Math.}, 63(9):1237--1247, 2010.

\bibitem[Bei97]{beig97}
R.~Beig.
\newblock T{T}-tensors and conformally flat structures on {$3$}-manifolds.
\newblock In {\em Mathematics of gravitation, {P}art {I} ({W}arsaw, 1996)},
  volume~41 of {\em Banach Center Publ.}, pages 109--118. Polish Acad. Sci.
  Inst. Math., Warsaw, 1997.

\bibitem[Cal19]{calegaririccinotes}
Danny Calegari.
\newblock Ricci flow (notes which are being transformed into chapter 6 of a
  book on 3-manifolds).
\newblock 2019.

\bibitem[CM05]{coldminiextinction}
Tobias~H. Colding and William~P. Minicozzi, II.
\newblock Estimates for the extinction time for the {R}icci flow on certain
  3-manifolds and a question of {P}erelman.
\newblock {\em J. Amer. Math. Soc.}, 18(3):561--569, 2005.

\bibitem[CMN]{cmn}
D.~Calegari, F.C. Marques, and A.~Neves.
\newblock Counting minimal surfaces in negatively curved 3-manifolds.
\newblock {\em arXiv:2002.01062 [math]}.

\bibitem[DeT83]{deturcktrick}
Dennis~M. DeTurck.
\newblock Deforming metrics in the direction of their {R}icci tensors.
\newblock {\em J. Differential Geom.}, 18(1):157--162, 1983.

\bibitem[Gut11]{guthlargeballs}
Larry Guth.
\newblock Volumes of balls in large {R}iemannian manifolds.
\newblock {\em Ann. of Math. (2)}, 173(1):51--76, 2011.

\bibitem[Ham82]{hamshorttime}
Richard~S. Hamilton.
\newblock Three-manifolds with positive {R}icci curvature.
\newblock {\em J. Differential Geometry}, 17(2):255--306, 1982.

\bibitem[JM87]{JM87}
Dennis Johnson and John~J. Millson.
\newblock Deformation spaces associated to compact hyperbolic manifolds.
\newblock In {\em Discrete groups in geometry and analysis ({N}ew {H}aven,
  {C}onn., 1984)}, volume~67 of {\em Progr. Math.}, pages 48--106.
  Birkh\"{a}user Boston, Boston, MA, 1987.

\bibitem[Kap94]{kap94}
Michael Kapovich.
\newblock Deformations of representations of discrete subgroups of {${\rm
  SO}(3,1)$}.
\newblock {\em Math. Ann.}, 299(2):341--354, 1994.

\bibitem[KM12a]{kmcounting}
Jeremy Kahn and Vladimir Markovi\'{c}.
\newblock Counting essential surfaces in a closed hyperbolic three-manifold.
\newblock {\em Geom. Topol.}, 16(1):601--624, 2012.

\bibitem[KM12b]{km}
Jeremy Kahn and Vladimir Markovic.
\newblock Immersing almost geodesic surfaces in a closed hyperbolic three
  manifold.
\newblock {\em Ann. of Math. (2)}, 175(3):1127--1190, 2012.

\bibitem[KY09]{knopfccf}
Dan Knopf and Andrea Young.
\newblock Asymptotic stability of the cross curvature flow at a hyperbolic
  metric.
\newblock {\em Proc. Amer. Math. Soc.}, 137(2):699--709, 2009.

\bibitem[Li]{chaolipolyhyp}
C.~Li.
\newblock Dihedral rigidity of parabolic polyhedrons in hyperbolic spaces.
\newblock {\em arXiv:2007.12563}.

\bibitem[MN12]{marquesnevesduke}
Fernando~C. Marques and Andr\'{e} Neves.
\newblock Rigidity of min-max minimal spheres in three-manifolds.
\newblock {\em Duke Math. J.}, 161(14):2725--2752, 2012.

\bibitem[Nun13]{nunes}
Ivaldo Nunes.
\newblock Rigidity of area-minimizing hyperbolic surfaces in three-manifolds.
\newblock {\em J. Geom. Anal.}, 23(3):1290--1302, 2013.

\bibitem[Per02]{perelman1}
G.~Perelman.
\newblock The entropy formula for the ricci flow and its geometric
  applications.
\newblock {\em arXiv:0211159 [math.DG]}, 2002.

\bibitem[Per03a]{perelman2}
G.~Perelman.
\newblock Finite extinction time for the solutions to the ricci flow on certain
  three-manifolds.
\newblock {\em arXiv:0307245 [math.DG]}, 2003.

\bibitem[Per03b]{perelman3}
G.~Perelman.
\newblock Ricci flow with surgery on three-manifolds.
\newblock {\em arXiv:0303109 [math.DG]}, 2003.

\bibitem[Sca02]{sca02}
Kevin~P. Scannell.
\newblock Local rigidity of hyperbolic 3-manifolds after {D}ehn surgery.
\newblock {\em Duke Math. J.}, 114(1):1--14, 2002.

\bibitem[Sep16]{seppi}
Andrea Seppi.
\newblock Minimal discs in hyperbolic space bounded by a quasicircle at
  infinity.
\newblock {\em Comment. Math. Helv.}, 91(4):807--839, 2016.

\bibitem[Son18]{songembness}
Antoine Song.
\newblock Embeddedness of least area minimal hypersurfaces.
\newblock {\em J. Differential Geom.}, 110(2):345--377, 2018.

\bibitem[SU82]{sacksuhl}
J.~Sacks and K.~Uhlenbeck.
\newblock Minimal immersions of closed {R}iemann surfaces.
\newblock {\em Trans. Amer. Math. Soc.}, 271(2):639--652, 1982.

\bibitem[SY79]{schoenyauincompressible}
R.~Schoen and Shing~Tung Yau.
\newblock Existence of incompressible minimal surfaces and the topology of
  three-dimensional manifolds with nonnegative scalar curvature.
\newblock {\em Ann. of Math. (2)}, 110(1):127--142, 1979.

\end{thebibliography}
	\bibliographystyle{alpha}
	
\end{document}